\documentclass[12pt]{article}
\usepackage{hyperref}
\usepackage{amsfonts}
\usepackage{enumerate}
\usepackage{graphicx}
\begin{document}

\begin{center}
{\large\bf A note on the relation between a singular linear discrete time system and a singular linear system of fractional nabla difference equations}

\vskip.20in

Charalambos P. Kontzalis$^{1}$ and\ Grigoris Kalogeropoulos$^{2}$\\[2mm]
{\footnotesize
$^{1}$Department of Informatics, Ionian University, Corfu, Greece\\[5pt]
$^{2}$Department of Mathematics, University of Athens, Greece}
\end{center}

{\footnotesize
\noindent
\textbf{Abstract:} In this article we focus our attention on the relation between a singular linear discrete time system and a singular linear system of fractional nabla difference equations whose coefficients are square constant matrices. By using matrix pencil theory, first we give necessary and sufficient condition to obtain a unique solution for the continuous time model. After by assuming that the input vector changes only at equally space sampling instants, we shall derive the corresponding discrete time state equation which yield the values of the solutions of the continuous time model which will connect the initial system to the singular linear system of fractional nabla difference equations.
.\\
\\[3pt]
\\[3pt]

\vskip.2in
\section{Introduction}

Many authors have studied generalized continuous \& discrete time systems, see [3-6, 10-12, 14-16, 19, 23, 24, 29-34, 37-39, 44-46], and their applications, see [2, 3, 20-22, 25, 26, 48-50, 53, 54]. Many of these results have already been extended to systems of differential \& difference equations with fractional operators, see [1, 13, 17, 18, 33-36, 40-42]. In this article, our purpose is to study the relation between a singular linear discrete time system and a singular linear system of fractional nabla difference equations whose coefficients are square constant matrices, into the mainstream of matrix pencil theory. We consider the singular matrix equation
\begin{equation}
FY'(t)=GY(t)+BV(t)
\end{equation}
with known initial conditions
\begin{equation}
Y(t_0)
\end{equation}        
where $F,G \in \mathcal{M}(m \times m;\mathcal{F})$, $B
\in \mathcal{M}(m \times r;\mathcal{F})$, $V(t)
\in \mathcal{M}(r \times 1;\mathcal{F})$, $Y(t)\in \mathcal{M}(m \times 1;\mathcal{F})$ (i.e., the algebra of
square matrices with elements in the field $\mathcal{F}$). For the
sake of simplicity, we set ${\mathcal{M}}_m = {\mathcal{M}}({m
\times m;\mathcal{F}})$ and ${\mathcal{M}}_{nm}  =
{\mathcal{M}}({n \times m;\mathcal{F}})$. We assume that the system (1) is singular, i.e. the matrix $F$ is singular.
\\\\
Difference equations of fractional order have recently proven to be valuable tools in the modeling of many phenomena in various fields of science and engineering. If we define $\mathbb{N}_\alpha$ by $\mathbb{N}_\alpha=\left\{\alpha, \alpha+1, \alpha+2,...\right\}$, $\alpha$ integer, and $n$ such that $0<n<1$ or $1<n<2$, then the nabla fractional operator in the case of Riemann-Liouville fractional difference of $n$-th order for any 
$Y_k:\mathbb{N}_a\rightarrow \mathbb{R}^{m\times1}$ is defined by, see [2-4, 29, 40],
\[
\nabla_\alpha^nY_k=\frac{1}{\Gamma(-n)}\sum^{k}_{j=\alpha}(k-j+1)^{\overline{-n-1}}Y_j,
\]
where the raising power function is defined by
\[
k^{\bar{\alpha}}=\frac{\Gamma(k+\alpha)}{\Gamma(k)}.
\]
We considered the following singular fractional discrete time system of the form
\begin{equation}
\begin{array}{cc}F\nabla_0^nY_k=GY_k+V_k, & k= 1, 2,...\end{array},
\end{equation}
Where  $F, G$ as defined for (1), $Y_k, V_k \in \mathcal{M}_{m1}({m \times
1;\mathcal{F}})$.
\\\\
\textbf{Definition 1.1.} Given $F,G\in \mathcal{M}_{rm}$ and an arbitrary $s\in\mathcal{F}$, the matrix
pencil $sF-G$ is called:
\begin{enumerate}
\item Regular when  $r=m$ and  det$(sF-G)\neq 0$;
\item Singular when  $r\neq m$ or  $r=m$ and det$(sF-G)\equiv 0$.
\end{enumerate}
In this article, we consider the case of the system (1) with a regular pencil. The class of $sF-G$ is characterized by a uniquely
defined element, known as complex Weierstrass canonical form, $sF_w -G_w$, see [23, 31, 33, 42], specified by the complete set of
invariants of $sF-G$. This is the set of elementary divisors obtained by
factorizing the invariant polynomials into powers of homogeneous polynomials irreducible over the field $\mathcal{F}$. In the case
where $sF-G$ is regular, we have elementary divisors of the following type:
\begin{itemize}
    \item elementary divisors of the type  $(s-a_j)^{p_j}$, \emph{are called finite elementary
    divisors}, where $a_j$ is a finite eigenvalue of algebraic multiplicity $p_j$;
    
    \item elementary divisors of the type  $\hat{s}^q=\frac{1}{s^q}$, are called \emph{infinite elementary divisors}, where $q$ is the algebraic multiplicity of the infinite eigenvalues.
\end{itemize}
We assume that $\sum_{i =1}^\nu p_j  = p$ and $p+q=m$.
\\\\
\textbf{Definition 1.2.} Let $B_1 ,B_2 ,\dots, B_l $ be elements of $\mathcal{M}_m$. The direct sum
of them denoted by $B_1  \oplus B_2  \oplus \dots \oplus B_l$ is
the blockdiag$\left[\begin{array}{cccc} B_1& B_2& \dots& B_l\end{array}\right]$.
\\\\ 
From the regularity of $sF-G$ there exist non-singular matrices $P$, $Q$ $\in \mathcal{M}_m$ such that 
\begin{equation}
PFQ = F_w  = I_p  \oplus H_q,
\]
\[
PGQ = G_w  = J_p  \oplus I_q.
\end{equation}
The complex Weierstrass form $sF_w -Q_w$ of the regular pencil $sF-G$ is defined by
\[
sF_w  - Q_w :=sI_p  - J_p  \oplus sH_q  - I_q,
\]
where the first normal Jordan-type element is uniquely defined by the set of the finite eigenvalues of $sF-G$ and has the form
\[
    sI_p  - J_p  := sI_{p_1 }  - J_{p_1 } (
    {a_1 }) \oplus  \dots  \oplus sI_{p_\nu  }  - J_{p_\nu  }.
    ({a_\nu  }).
\]
The second uniquely defined block $sH_q -I_q$ corresponds to the infinite eigenvalues of $sF-G$  and has the form
\[
    sH_q  - I_q  := sH_{q_1 }  - I_{q_1 }  \oplus 
    \dots  \oplus sH_{q_\sigma  }  - I_{q_\sigma}.
\]
The matrix $H_q$ is a nilpotent element of $\mathcal{M}_q$  with index
$q^* = \max \{ {q_j :j = 1,2, \ldots ,\sigma } \}$, i.e.
\[
    H^{q^*}_q=0_{q, q}.
\]
The matrices $I_{p_j } ,J_{p_j } ({a_j }),H_{q_j }$ are defined as
\[
   I_{p_j }  = \left[\begin{array}{ccccc} 
   1&0& \ldots & 0&0\\
   0& 1 &  \ldots&0 &0 \\
   \vdots & \vdots & \ddots & \vdots &\vdots \\
   0 & 0 & \ldots  & 0 &1
   \end{array}\right]
   \in {\mathcal{M}}_{p_j }, 
   \]
   \[
   J_{p_j } ({a_j }) =  \left[\begin{array}{ccccc}
   a_j  & 1 & \dots&0  & 0  \\
   0 & a_j  &   \dots&0  & 0  \\
    \vdots  &  \vdots  &  \ddots  &  \vdots  &  \vdots   \\
   0 & 0 &  \ldots& a_j& 1\\
   0 & 0 & \ldots& 0& a_j
   \end{array}\right] \in {\mathcal{M}}_{p_j },
\]
\[
 H_{q_j }  = \left[
\begin{array}{ccccc} 0&1&\ldots&0&0\\0&0&\ldots&0&0\\\vdots&\vdots&\ddots&\vdots&\vdots\\0&0&\ldots&0&1\\0&0&\ldots&0&0
\end{array}
\right] \in {\mathcal{M}}_{q_j }.
  \]
For algorithms about the computations of the Jordan matrices, see [28, 39, 39, 47].

\section{The solution of a singular linear state space equation}

Firstly, we will derive the corresponding discrete time state equation which yield the values of $Y(t)$ at $t_k= k\cdot T$ for $k\geq 0$. We will use the same method as used in [39]. In order to do so, we obtain formulas for the solutions of (1) with a regular matrix pencil and we
give necessary and sufficient conditions for existence and uniqueness of solutions.
\\\\
\textbf{Theorem 2.1.} Consider the system (1), (2) and let the $p$ linear independent (generalized) eigenvectors of the finite eigenvalues of the pencil $sF-G$ be the columns of a matrix $Q_p$. Then the solution is unique if and only if
\begin{equation}
Y(t_0)\in colspan Q_p+QK(t_0)
\end{equation}
Moreover the analytic solution is given by
\begin{equation}
    Y(t)=Q_pe^{J_p(t-t_0)}Z_p(t_0)+QK(t)
\end{equation}
where $K(t)=\left[\begin{array}{c} \int^{t}_{t_0}e^{J_p(t-s)}B_pV(s)ds \\ -\sum^{q_{*}-1}_{i=0}H_q^iB_qV^i(t)\end{array}\right]$ and $PB=\left[\begin{array}{c} B_p \\ B_q\end{array}\right]$, with $B_p \in \textsl{M}_{pm}$, $B_q \in \textsl{M}_{qm}$.
\\\\
\textbf{Proof.} Consider the transformation
\begin{equation}
    Y(t)=QZ(t)
\end{equation}
Substituting the previous expression in to (1) we obtain
\[
    FQZ'(t)=GQZ(t)+BV(t).
\]
Whereby, multiplying by $P$, we arrive at
\[
    F_wZ(t)=G_w Z(t)+PBV(t).
\]
Moreover, we can write $Z(t)$ as
$Z(t) =\left[\begin{array}{c}
 Z_p (t) \\
 Z_q (t)
 \end{array}\right] $.
Taking into account the above expressions, we arrive easily at two subsystems of (1). The subsystem
\begin{equation}
    Z_p'(t) = J_p Z_p (t)+B_pV(t) 
 \end{equation}
and the subsystem
\begin{equation}
    H_q Z_q'(t) = Z_q (t)+B_q V(t)
\end{equation}
The subsystem (8) has the unique solution 
\begin{equation}
    Z_p(t)=e^{J_p(t-t_0)}Z_p(t_0)+\int^{t}_{t_0}e^{J_p(t-s)}B_pV(s)ds, t\geq t_0.
\end{equation}
see [2, 7, 8-10, 48-52]. By applying the Laplace transform we get the solution of subsystem (9)  
\begin{equation}
Z_q(t)=-\sum^{q_*-1}_{i=0}H_q^iB_qV^{(i)}(t)
\end{equation}
 Let $Q=\left[\begin{array}{cc} Q_p& Q_q\end{array} \right]$, where $Q_p \in \textsl{M}_{mp}$, $Q_q\in \textsl{M}_{mq}$ the matrices with columns the $p$, $q$ generalized eigenvectors of the finite and infinite eigenvalues respectively. Then we obtain
\[
     Y(t) = QZ(t) =
     [Q_p  Q_q ]
     \left[\begin{array}{c}
     e^{J_p(t-t_0)}Z_p(t_0)+\int^{t}_{t_0}e^{J_p(t-s)}B_pV(s)ds  \\
     -\sum^{q_{*}-1}_{i=0}H_q^iB_qV^{(i)}(t),
    \end{array}\right]  
    \]
    \[
    Y(t)=Q_pe^{J_p(t-t_0)}Z_p(t_0)+\int^{t}_{t_0}e^{J_p(t-s)}B_pV(s)ds-Q_q\sum^{q_{*}-1}_{i=0}H_q^iB_qV^{(i)}(t).
    \]
    \[
    Y(t)=Q_pe^{J_p(t-t_0)}Z_p(t_0)+QK(t).
    \]
The solution that exists if and only if 
\[
Y(t_0)=Q_pZ_p(t_0)+QK(t_0),
\]
or
\[
Y(t_0)\in colspan Q_p+QK(t_0).
\]
\textbf{Proposition 2.1.} Assume the system (1) with initial conditions (2). Then if (5) holds the unique solution is given from the formula
\begin{equation}
Y(t)=F(t,t_0)Y(t_0)+\int^{t}_{t_0}F(t,s)Q_pB_pV(s)ds+Q_q\sum^{q_{*}-1}_{i=0}H_q^iB_q(V^i(t_0)-V^i(t)),
\end{equation}
where
\begin{equation}
F(t,t_0)=Q\left[\begin{array}{cc}e^{J_p(t-t_0)} & 0_{p, q} \\ 0_{q, p} & I_p \end{array} \right]Q^{-1}
\end{equation}
and
\begin{equation}
F(t,s)=Q\left[\begin{array}{cc}e^{J_p(t-s)} & 0_{p, q} \\ 0_{q, p} & I_p \end{array} \right]Q^{-1}.
\end{equation}
\textbf{Proof.} From (6) the solution of (1) is 
\[
Y(t)=Q_pe^{J_p(t-t_0)}Z_p(t_0)+\int^{t}_{t_0}Q_pe^{J_p(t-s)}B_pV(s)ds-Q_q\sum^{q_{*}-1}_{i=0}H_q^iB_qV^{(i)}(t),
\]
or,
\[
Y(t)=Q_pe^{J_p(t-t_0)}Z_p(t_0)+Q_qZ_q(t_0)-Q_qZ_q(t_0)+
\]
\[
+\int^{t}_{t_0}Q_pe^{J_p(t-s)}B_pV(s)ds-Q_q\sum^{q_{*}-1}_{i=0}H_q^iB_qV^{(i)}(t),
\]
or
\begin{equation}
Y(t)=\left[\begin{array}{cc}Q_p&Q_q\end{array}\right]\left[\begin{array}{cc}e^{J_p(t-t_0)}&0_{p, q}\\0_{q, p}&I_q\end{array}\right]\left[\begin{array}{c}Z_p(t_0)\\Z_q(t_0)\end{array}\right]+
\]
\[
+\int^{t}_{t_0}Q_pe^{J_p(t-s)}B_pV(s)ds+Q_q(-Z_q(t_0)-\sum^{q_{*}-1}_{i=0}H_q^iB_qV^{(i)}(t)).
\end{equation}   
By substituting $\left[\begin{array}{c}Z_p(t_0)\\Z_q(t_0)\end{array}\right]=Q^{-1}Y(t_0)$ and $-Z_q(t_0)=\sum^{q_{*}-1}_{i=0}H_q^iB_qV^{(i)}(t_0)$, we get
\begin{equation}
Y(t)=F(t,t_0)Y(t_0)+\int^{t}_{t_0}Q_pe^{J_p(t-s)}B_pV(s)ds+Q_q\sum^{q_{*}-1}_{i=0}H_q^iB_q(V^{(i)}(t_0)-V^{(i)}(t))   
\end{equation}
Let $Q^{-1}$ be the inverse matrix of $Q$ defined as $Q^{-1}=\left[\begin{array}{c}\bar Q_p\\\bar Q_q\end{array}\right]$, with $\bar Q_p \in \textsl{M}_{pn}$ and $\bar Q_q \in \textsl{M}_{qn}$. Then
\[
Q^{-1}Q=\left[\begin{array}{c}\bar Q_p\\\bar Q_q\end{array}\right]\left[\begin{array}{cc}Q_p& Q_q\end{array}\right]=\left[\begin{array}{cc}\bar Q_pQ_p&\bar Q_pQ_q\\\bar Q_qQ_p&\bar Q_qQ_q\end{array}\right].
\]
From where we get $\bar Q_pQ_p=I_p$, $\bar Q_pQ_q=0_{p, q}$, $\bar Q_qQ_p=0_{p, q}$, $\bar Q_qQ_p=I_q$. Then
\[
F(t,s)Q_p=Q\left[\begin{array}{cc}e^{J_p(t-t_0)}&0_{p, q}\\0_{q, p}&I_q\end{array}\right]Q^{-1}Q_p,
\]
or
\[
F(t,s)Q_p=\left[\begin{array}{cc}Q_p&Q_q\end{array}\right]\left[\begin{array}{cc}e^{J_p(t-t_0)}&0_{p, q}\\0_{q, p}&I_q\end{array}\right]\left[\begin{array}{c}I_p\\0_{q, p}\end{array}\right],
\]
or
\[
F(t,s)Q_p=Q_pe^{J_p(t-s)}
\]
and thus
\[.
Y(t)=F(t,t_0)Y(t_0)+\int^{t}_{t_0}F(t,s)Q_pB_pV(s)ds+Q_q\sum^{q_{*}-1}_{i=0}H_q^iB_q(V^i(t_0)-V^i(t)).
\]
In the following analysis, we will use the notation $V(t)=V(kT)$ and $Y_k=Y(kT)$. We assume that the input $V(t)$ is sampled and fed to a zero order hold, so that all the components of $V(t)$ are constant over the interval between any two consecutive sampling instants, or for $kT\leq t \leq (k+1)T$. Thus $V(t)$ changes only during the time periods $\begin{array}{cc} t_k=kT, & k=0, 1, 2,... \end{array}$. For $t=kT$, (12) takes the form
\begin{equation}
Y_k=F(kT,0)Y_0+\int^{kT}_{0}F(kT,s)Q_pB_pV(s)ds+Q_q\sum^{q_{*}-1}_{i=0}H_q^iB_q(V^{(i)}(0)-V^{(i)}(kT))
\end{equation}
and for $t=(k+1)T$ in (12),
\[
Y_{k+1}=F((k+1)T,0)Y_0+\int^{(k+1)T}_{0}F((k+1)T,s)Q_pB_pV(s)ds+
\]
\[
+Q_q\sum^{q_{*}-1}_{i=0}H_q^iB_q(V^{(i)}(0)-V^{(i)}((k+1)T)).
\]
By multiplying (17) with $F(T,0)$ and taking into account that 
\[
F(T,0)F(kT,0)=F((k+1)T,0),\quad F(T,0)F(kT,s)=F((k+1)T,s)
\]
and
\[
F(T,0)Y_k=F((k+1)T,0)Y_0+\int^{kT}_{0}F((k+1)T,s)Q_pB_pV(s)ds+
\]
\[
+F(T,0)Q_q\sum^{q_{*}-1}_{i=0}H_q^iB_q(V^{(i)}(0)-V^{(i)}(kT)).
\]
In addition,
\[
Y_{k+1}-F(T,0)Y_k=\int^{(k+1)T}_{kT}F((k+1)T,s)Q_pB_pV(s)ds+
\]
\[
+Q_q\sum^{q_{*}-1}_{i=0}H_q^iB_q(V^{(i)}(0)-V^{(i)}((k+1)T))-F(T,0)Q_q\sum^{q_{*}-1}_{i=0}H_q^iB_q(V^{(i)}(0)-V^{(i)}(kT)).
\]
We should note that $V(t)$ may jump at $t=(k+1)T$ and thus $V((k+1)T)$ may be different from $V(kT)$. Such a jump in $V(t)$ at $t=(k+1)T$, the upper limit of integration, does not affect the value of the integral in this last equation, because the integrand does not involve impulse functions. Thus by considering the substitution $s=kT+w$ in the following integral 
\[
\int^{(k+1)T}_{kT}F((k+1)T,s)Q_pB_pV(s)ds=\int^{(k+1)T}_{0}F((k+1)T,kT+w)Q_pB_pV(kT+w)dw,
\]
or
\[
\int^{(k+1)T}_{kT}F((k+1)T,s)Q_pB_pV(s)ds=\int^{(k+1)T}_{0}F((k+1)T,kT+w)dwQ_pB_pV_k
\]
and by using the relations
\[
\int^{(k+1)T}_{kT}F((k+1)T,s)Q_pB_pV(s)ds=\int^{(k+1)T}_{0}F((k+1)T,kT)F(0,w)dwQ_pB_pV_k,
\]
or
\[
\int^{(k+1)T}_{kT}F((k+1)T,s)Q_pB_pV(s)ds=\int^{(k+1)T}_{0}F(T,0)F(kT,kT)F(0,w)dwQ_pB_pV_k,
\]
and since $F(kT,kT)=I_m$
\[
\int^{(k+1)T}_{kT}F((k+1)T,s)Q_pB_pV(s)ds=\int^{(k+1)T}_{0}F(T,0)F(0,w)dwQ_pB_pV_k.
\]
Furthermore since $F(T,0)F(0,w)=F(T,w)=F(T-w,0)=F(l,0)$, for $l=T-w$,
\[
\int^{(k+1)T}_{kT}F((k+1)T,s)Q_pB_pV(s)ds=\int^{(k+1)T}_{0}F(l,0)dwQ_pB_pV_k.
\]
Using the approximation $V_k^{(i)}=V^{(i)}(kT)\cong\frac{\sum^{i}_{j=0}(-1)^j(\begin{array}{c}i\\j\end{array})V_{k-i+j}}{T^i}$
and the obvious relation $V^{(i)}(0)-V^{(i)}(kT)=(-T)V^{(i+1)}_{k+1}$, by setting $
V_{k+1}^{(i+1)}\cong\frac{\sum^{i+1}_{j=0}(-1)^j(\begin{array}{c}i+1\\j\end{array})V_{k-i+j}}{T^{i+1}}$ we have
\[
V_k^{(i)}-V_{k+1}^{(i+1)}\cong\frac{\sum^{i+1}_{j=0}(-1)^{j+1}(\begin{array}{c}i+1\\j\end{array})V_{k-i+j}}{T^i}.
\]
Then 
\begin{equation}
Y_{k+1}\cong F(T,0)Y_k+\int^{(k+1)T}_{0}F(l,0)dwQ_pB_pV_k+Q_q\sum^{q_{*}-1}_{i=0}H_q^iB_q\frac{\sum^{i+1}_{j=0}(-1)^{j+1}(\begin{array}{c}i+1\\j\end{array})V_{k-i+j}}{T^i}.
\end{equation}
\textbf{Remark 2.1.} The term $V_{k+i}$ in the state of discrete time system is just in correspondence with the term $V^{(i)}(t)$ in continuous time system.
\\\\
Let 
\[
A_k=F(T,0),\quad U_k=\int^{(k+1)T}_{0}F(l,0)dwQ_pB_pV_k+Q_q\sum^{q_{*}-1}_{i=0}H_q^iB_q\frac{\sum^{i+1}_{j=0}(-1)^{j+1}(\begin{array}{c}i+1\\j\end{array})V_{k-i+j}}{T^i}
\]
and
\[
Y_{k+1}=A_kY_k+U_k.
\]
Furthermore,
\[
\begin{array}{c}
     Y_k=A_{k-1}Y_{k-1}+U_{k-1},\\
     Y_{k-1}=A_{k-2}Y_{k-2}+U_{k-2},\\
     Y_{k-2}=A_{k-3}Y_{k-3}+U_{k-3},\\
		Y_{k-4}=A_{k-4}Y_{k-4}+U_{k-4},\\
    \vdots\\
     Y_2=A_1Y_1+U_1,\\\\
		Y_1=A_0Y_0+U_0.
    \end{array}
   \]
By taking the sum of the above equations, 
\[
Y_k=(A_{k-1}-I_m)Y_{k-1}+(A_{k-2}-I_m)Y_{k-2}+...+(A_1-I_m)Y_1+A_0Y_0+U_{k-1}+...+U_1+U_0,
\]
for 
\[
A_{k-j}-I_m=\frac{1}{\Gamma(-n)}(k-j)^{\overline{-n-1}}F,\quad V_k=U_{k-1}+...+U_1+U_0,
\]
we arrive at (3).

\section*{Conclusions}
In this article we focused our attention on the relation between a singular linear discrete time system and a singular linear system of fractional nabla difference equations whose coefficients are square constant matrices. We derived the corresponding discrete time matrix equation which yield the values of the solutions of the continuous time model which connects the initial system to a singular linear system of fractional nabla difference equations.

\end{document}